\documentclass{llncs} 
\usepackage{amsmath,epsf}
\usepackage[latin1]{inputenc}
\usepackage{amsfonts}

\title{Bases explicites et conjecture $n!$}
\author{Jean-Christophe Aval}
\institute{Laboratoire A2X, Universit\'e Bordeaux 1\\ 351 cours de la Lib\'eration, F-33405 Talence cedex\\ e-mail : {\tt aval@math.u-bordeaux.fr}}

\begin{document}

\maketitle

\begin{abstract}
Le but de ce travail est d' obtenir pour l'espace $M_{\mu}$, relatif la conjecture $n!$, une base explicite et monomiale. Ce but est atteint dans le cas des partitions qui ont la forme d'une équerre, i.e. $\mu=(K+1,1^L)$. Nous introduisons en effet une famille pour laquelle nous démontrons qu'elle est de cardinal $n!$, qu'elle est libre et qu'elle engendre $M_{\mu}$. Nous déduisons de cette étude une base simple pour $I_{\mu}$, l'idéal annulateur de $\Delta_{\mu}$. Cette méthode permet aussi de donner de façon directe une base monomiale pour le sous-espace de $M_{\mu}$ formé des éléments de degré en $x$ nul. \\
{\bf Abstract} 
The aim of this work is to construct a monomial and explicit basis for the space $M_{\mu}$ relative to the $n!$ conjecture. We succeed completely for hook-shaped partitions, i.e. $\mu=(K+1,1^L)$. We are indeed able to exhibit a basis and to verify that its cardinality is $n!$, that it is linearly independent and that it spans $M_{\mu}$. We deduce from this study an explicit and simple basis for $I_{\mu}$, the annulator ideal of $\Delta_{\mu}$. This method is also successful for giving directly a basis for the homogeneous subspace of $M_{\mu}$ consisting of elements of $0$ $x$-degree.
\end{abstract}

\section{Introduction}

Soit $\mu=(\mu_1\ge\mu_2\ge\dots\ge\mu_k>0)$ une partition de l'entier $n$. Nous identifierons $\mu$ et son diagramme de Ferrers (en utilisant la convention ``Française''). \`A chaque cellule $s$ du diagramme de Ferrers, nous associons ses coordonnées $(i,j)$, o\`u $i$ est la hauteur de $s$ et $j$ la position de $s$ dans sa ligne. Les paires $(i-1,j-1)$ apparaissant lorsque $s$ décrit toutes les cellules de $\mu$ seront désignées comme les biexposants de $\mu$. Soient désormais $(p_1,q_1),\ldots,(p_n,q_n)$ ces biexposants rangés selon l'ordre lexicographique et introduisons :
$$\Delta_{\mu}(x,y)=\Delta_{\mu}(x_1,\ldots,x_n;y_1,\ldots,y_n)=\det(x_i^{p_j}y_i^{q_j})_{i,j=1\ldots n}.$$
Soit $M_{\mu}$ l'ensemble des polynômes en $x_1,\ldots,$ $x_n;$ $y_1,\ldots,y_n$ qui peuvent s'écrire comme combinaison linéaire des dérivées partielles de $\Delta_{\mu}$, ce que l'on peut écrire sous la forme : 
$$M_{\mu}={\cal L}\{\partial_x^p\partial_y^q\Delta_{\mu}(x,y)\}$$
o\`u $\partial_x^p=\partial_{x_1}^{p_1}\ldots\partial_{x_n}^{p_n}$ et $\partial_y^q=\partial_{y_1}^{q_1}\ldots\partial_{y_n}^{q_n}$. On peut alors énoncer la conjecture $n!$ ainsi :

\noindent
{\bf Conjecture 1 (conjecture $n!$).} {\it Soit $\mu$ une partition de $n$, alors on a :\\ $\dim M_{\mu}=n!$}.

Cette conjecture a été énoncée pour la première fois par A. Garsia et M. Haiman. Elle est centrale dans le cadre de leur étude des polynômes de Macdonald (cf. [6], [7]). Plus précisément, Macdonald a introduit dans [13] une nouvelle base de l'espace des fonctions symétriques, et des coefficients associés, appelés polynômes de Macdonald-Kostka, $K_{\lambda\mu}(q,t)$, qui sont a priori des fractions rationnelles en $q,t$. Macdonald a conjecturé que :
  
\noindent
{\bf Conjecture 2 (conjecture MPK).} {\it Les fonctions $K_{\lambda\mu}(q,t)$ sont des polynômes à coefficients entiers positifs.}

A. Garsia et M. Haiman, en cherchant une interprétation de ces fonctions $K_{\lambda\mu}(q,t)$ liée à la théorie des représentations, firent la conjecture suivante :

\noindent
{\bf Conjecture 3 (conjecture $C=\tilde H$).} {\it Pour l'action diagonale de $S_n$, $M_{\mu}$ est une version de la représentation régulière à gauche. De plus, si on note $C_{\lambda\mu}(q,t)$ la multiplicité bigraduée du caractère $\chi_{\lambda}$ dans le caractère bigradué de $M_{\mu}$ alors : $C_{\lambda\mu}(q,t)=K_{\lambda\mu}(q,1/t)t^{n(\mu)}$ o\`u $n(\mu)=\sum_{i=1}^k(i-1)\mu_i$.}  

La Conjecture 3 implique clairement les Conjectures 1 et 2. Il est tout à fait remarquable que M. Haiman a récemment montré, en utilisant la théorie des schémas de Hilbert que la conjecture $n!$ implique la conjecture $C=\tilde H$. Une partie de la conjecture MPK est le fait que les fonctions $K_{\lambda\mu}$ sont des polynômes, ce qui n'est pas évident vu leur définition. Cette partie a été prouvée récemment dans plusieurs articles indépendants (cf. [8], [9], [11], [12], [15]).

De plus, M. Haiman a très récemment obtenu une preuve de la conjecture $n!$, impliquant ainsi la conjecture $C=\tilde H$. Cependant la recherche d'une preuve combinatoire et l'obtention de bases explicites maintiennent à ces questions un intérêt certain.

Quand $\mu=(1^n)$ ou $\mu=(n)$, $\Delta_{\mu}$ est simplement le déterminant de Vandermonde en $x$ et $y$ respectivement. Dans ce cas, c'est un résultat classique (cf. [3]) que $\dim M_{\mu}=n!$ et on obtient de plus dans ce cas une base explicite simple pour $M_{\mu}$. Des preuves combinatoires de plusieurs cas particuliers se trouvent dans [1], [5], [7], [14]. 

Dans ce travail, notre but est de proposer une nouvelle méthode visant à démontrer de façon combinatoire la conjecture $n!$ dans certains cas particuliers et d'obtenir des bases ``simples''. Nous cherchons à construire une base explicite pour $M_{\mu}$, constituée de dérivées monomiales de $\Delta_{\mu}$. Nous présentons ici comment nous sommes capables de le faire pour les équerres, i.e. $\mu=(K+1,1^L)$ avec $K+L+1=n$. Dans la deuxième partie nous décrivons la façon de construire cette famille et prouvons que son cardinal est bien $n!$. Dans la troisième partie, nous démontrons que notre famille engendre $M_{\mu}$. De plus, nous déduisons de cette preuve une base explicite et simple pour $I_{\mu}$, l'idéal annulateur de $\Delta_{\mu}$. Dans la quatrième partie nous montrons, et ce de façon totalement nouvelle, que les éléments de notre famille sont linéairement indépendants. Dans la cinquième et denière partie nous expliquons comment cette méthode est également efficace pour le sous-espace de $M_{\mu}$ constitué des éléments de degré en $x$ nul. Nous obtenons en fait un procédé construisant directement une base pour cet espace. 

Cet article est une version préliminaire de l'article [2]. En particulier, les preuves ne sont ici qu'esquissées. Le but est ici de donner un exposé des résultats et des méthodes développés dans cet article. Le détail complet des démonstra\-tions se trouve dans [2].
 
\section{Construction et énumeration}

Soit $\mu$ une partition de $n$ dont le diagramme de Ferrers est une équerre, i.e. $\mu=(K+1,1^L)$ avec $K+L+1=n$.

\subsection{Construction}

Donnons-nous un axe. Un dessin associé à $\mu$ est un ensemble de $K$ colonnes situées au-dessus de l'axe (appelées colonnes-$y$) de hauteur $K, K-1,\ldots, 1$, rangées en taille décroissante de la gauche vers la droite, et de $L$ colonnes, éventuellement intercalées, au-dessous de l'axe (appelées colonnes-$x$) de profondeur $L, L-1,\ldots , 1$, rangées par profondeur décroissante de la gauche vers la droite.

Voici un exemple de dessin :
\vskip 0.5 cm

\centerline{
\epsffile{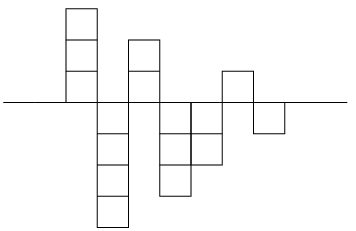}}

\vskip 0.3 cm

associé à la partition:

\centerline{
\epsffile{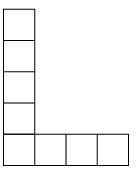}}

\vskip 0.3 cm

Nous allons maintenant mettre des croix dans les cases des dessins. Comme nous ne distinguerons pas deux dessins ayant le même nombre de croix dans chaque colonne, nous plaçons les croix près de l'axe. Les règles pour placer des croix dans un dessin sont les suivantes : 
\begin{enumerate} 
\item le nombre de croix dans les colonnes-$x$ est quelconque (limité seulement par la taille de la colonne) ;
\item le nombre de croix dans les colonnes-$y$ dépend des croix-$x$. Pour une colonne n'ayant à sa droite aucune colonne-x, le nombre de croix est arbitraire (limité seulement par la hauteur de la colonne). Dans l'autre cas, on regarde la première colonne-$x$ à droite de notre colonne-$y$ qui est ``unie'' (i.e. qui ne contient que des croix ou des cases blanches). Il y en a toujours une, au moins celle de profondeur 1. Alors :
\begin{itemize}
\item si elle est toute blanche, on impose au moins une croix dans la colonne-$y$ ;
\item si elle est pleine de croix, on impose au moins un blanc dans la colonne-$y$.
\end{itemize}
\end{enumerate}

\noindent {\it Remarque 1.} La famille de dessins que nous avons définie est invariante par interversion des cases croisées et des cases blanches. On appelle flip l'opérateur correspondant (il est different du flip introduit par A. Garsia et M. Haiman dans [7], que nous notons désormais Flip).

Donnons un exemple de dessin avec croix pour la même partition que ci-dessus :
\vskip 0.5 cm

\centerline{
\epsffile{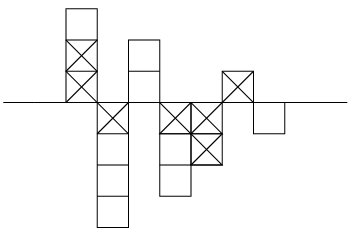}}

\vskip 0.3 cm
 
On associe maintenant aux dessins (avec croix) des opérateurs de dériva\-tion. On donne un indice aux places sur le dessin de la gauche vers la droite et de $1$ à $n-1$. Alors pour chaque croix-$x$ en place $i$, on dérive une fois par rapport à $x_i$, et on procède même pour les $y$. Par exemple pour le dernier dessin l'opérateur de dérivation associé est : $\partial_D= \partial y_1^2\partial x_2 \partial x_4\partial x_5^2\partial y_6$.

\subsection{\'Enumeration}

Nous notons ${\cal D}$ l'ensemble des dessins définis au paragraphe précédent et nous allons commencer par vérifier que son cardinal est $n!$.

Comme le nombre de choix ne dépend pas des croix-$x$ mais seulement de la forme du dessin et plus précisément du nombre de colonnes-$y$ à droite de la colonne-$x$ de profondeur 1 (notons $k_1$ ce nombre), nous pouvons écrire que le cardinal est :
$$\sum_{k_1+k_2=K} 2.3\ldots(k_1+1).(k_1+1)\ldots(k_1+k_2).(L+1)!\ {{k_2+L-1} \choose {k_2}}$$
$$=L(L+1)K! \sum_{k_2=0}^K \frac {(k_2+L-1)!} {k_2!} (K+1-k_2)$$
$$=(L+1)!K!\sum_{k_2=0}^K {{L-1+k_2}\choose{L-1}} {{K+1-k_2}\choose{1}}$$
$$=(L+1)!K!{{K+L+1}\choose{L+1}}=(K+L+1)!$$
d'après la formule de Chu-Vandermonde ([4], p. 163).

\section{Preuve que la famille engendre $M_{\mu}$}

Nous montrons ici que $\{\partial_D \Delta_{\mu}\}_{D\in{\cal D}}$ engendre $M_{\mu}$, et pour ce faire, nous commençons par étudier $I_{\mu}$, l'idéal annulateur de $\Delta_{\mu}$.

\subsection{\'Etude de $I_{\mu}$}

Pour un polynôme $P$, nous écrivons $P\equiv 0$ si $P(\partial)\Delta_{\mu}=0$, i.e. $P\in I_{\mu}$. Nous notons $h_k$ la $k$-ième fonction symétrique homogène complète. Notons aussi $X$ une partie de $(x_1,x_2,\ldots,x_n)$, $Y$ une partie de $(y_1,y_2,\ldots,y_n)$, $|X|$ et $|Y|$ leurs cardinaux. Posons également $\bar{X}=\prod_{x\in X}x$ et $\bar{Y}=\prod_{y\in Y}y$.

Il est aisé de constater que :
\begin{enumerate}
\item pour tout $1\le i\le n$, $x_iy_i\equiv 0$;
\item $\bar{X}\equiv 0$ dès que $|X|>L$;
\item $\bar{Y}\equiv 0$ dès que $|Y|>K$;
\item pour tout polynôme $P$ symétrique homogène non constant, $P\equiv 0$.
\end{enumerate}

\noindent Grâce à ces observations, nous pouvons établir les propositions suivantes :

\vskip 0.2 cm
\noindent {\bf Proposition 1.} 
{\it $h_{k}(Y)\equiv 0$
dès que $k>0$ et $k+|Y|>n$.}

\vskip 0.2 cm
\noindent {\bf Proposition 2.} 
{\it $\bar{Y}h_{k}(Y')\equiv 0$
dès que $k>0$, $k+|Y|>K$ et $Y\subset Y'$.}

\vskip 0.2 cm
\noindent {\bf Proposition 3.} 
{\it $h_{k}(Y)h_{l}(X)\equiv 0$
dès que $k>0$, $l>0$ et 
\begin{itemize}
\item soit $Y\subset X$ et $k+l+|Y|>n$,
\item soit $X\subset Y$ et $k+l+|X|>n$.
\end{itemize}
}

Toutes ces propositions se démontrent à partir des observations 1-2-3-4 par de simples récurrences.

\subsection{Application}

Nous utilisons les propositions précédentes pour montrer que toute dérivée monomiale de $\Delta_{\mu}$ est une combinaison linéaire des dérivées : $\{\partial_D\Delta_{\mu}\}_{D \in {\cal D}}$. 

\vskip 0.3 cm

\noindent {\bf Théorème 1. }{\it $\{\partial_D\Delta_{\mu}\}_{D \in {\cal D}}$ engendre $M_{\mu}$.}

\vskip 0.3 cm

Nous montrons en fait que toute dérivée monomiale qui n'est pas associé à un dessin de ${\cal D}$ peut s'écrire comme combinaison linéaire de monômes strictement plus petits pour l'ordre lexicographique sur $x_1,\dots,x_n,y_1,\dots,y_n$. Pratiquement, pour tout monôme que l'on souhaite éliminer, les propositions précédentes d'écrire un polynôme de $I_{\mu}$ dont le terme dominant pour l'ordre lexicographique est le monôme en question.

\subsection{Conclusion}

Nous pouvons déduire de ce qui précède et du fait que notre famille est libre (ce qui est l'objet du prochain paragraphe) que les polynômes décrits aux points 1-2-3-4 au début de l'étude de $I_{\mu}$ forment une base de celui-ci :

\vskip 0.3 cm
\noindent {\bf Théorème 2.} 
{\it Notons $<G>$ l'idéal engendré par un ensemble $G$, alors pour $\mu=(K+1,1^L)$, partition de $n$, nous avons :
$$I_{\mu}=<h_i(X_n),\ 1\le i\le n;\ h_i(Y_n),\ 1\le i\le n;$$ 
$$x_iy_i,\ 1\le i\le n;\ X,\ |X|=L+1;\ Y,\ |Y|=K+1>.$$}

\section{Preuve de l'indépendance}

\subsection{Exposition et réduction du problème}

Nous voulons montrer que notre famille est linéairement indépendante. 

Pour un dessin $D$ nous notons $S$ l'ensemble des cases croisées et $T$ celui des cases blanches. Nous associons à ces deux objets leurs opérateurs de dérivation $\partial_S$ et $\partial_T$, de façon évidente. Nous notons de plus ${\cal S}$ l'ensemble des $S$ pour $D$ dans ${\cal D}$. Nous vérifions assez facilement que deux dessins distincts donnent des $S$ et des $T$ différents. Nous sommes même capables de reconstruire $D$ à partir de $S$ ou $T$.

\vskip 0.3 cm
\noindent {\bf Théorème 3.} 
{\it La famille $\{\partial_S.\Delta_{\mu}\}_{S \in {\cal S}}$ est linéairement indépendante.}
 
\vskip 0.2 cm
Pour ce faire, nous introduisons quelques définitions. Soient $D=(S,T)$ et $D_1=(S_1,T_1)$ deux dessins différents. Nous dirons que $D_1$ est un fils de $D$ si $\partial_T\circ\partial_{S_1}.\Delta_{\mu}\in {\mathbb Z}^*$. Nous noterons $T+S_1$ le résultat de la superposition place par place des cases de $T$ et de $S_1$. En itérant la filiation, nous obtenons la notion de descendant. Ceci étant posé, on obtient facilement le :

\vskip 0.3 cm
\noindent {\bf Lemme 1.}
{\it Pour prouver l'indépendance, il est suffisant de montrer qu'un dessin ne peut être son propre descendant (i.e. il n'y a pas de ``boucles'').}

\subsection{Définition de la complétude}

Soit $D_1$ un dessin et $D_2$ l'un de ses fils. Nous dirons qu'il y a complétude sur les $k$ premières places si les hauteurs des colonnes-$y$ de $T_1+S_2$ lues de la gauche vers la droite sont $K,K-1,K-2,\ldots$ et si on a la même chose pour les colonnes-$x$.

Afin d'obtenir une caractérisation plus quantitative de la complétude, nous regardons les parties gauches des dessins (i.e. les $k-1$ premières places). Nous définissons $d$ comme la différence entre le nombre de fois o\`u une colonne-$y$ de $D_1$ a été remplacée par une colonne-$x$ blanche dans $D_2$ et le nombre de fois o\`u une colonne-$x$ de $D_1$ a été remplacée dans $D_2$ par une colonne-$y$ blanche. Nous définissons de même $d'$ comme la différence entre le nombre de fois o\`u une colonne-$y$ croisée de $D_1$ a été remplacée par une colonne-$x$ dans $D_2$ et le nombre de fois o\`u une colonne-$x$ croisée de $D_1$ a été remplacée dans $D_2$ par une colonne-$y$. Notons que $d$ et $d'$ sont relatifs à $k-1$. 

Nous introduisons aussi les notations suivantes : $b_1$ (resp. $b_2$) représente le nombre de cases blanches en place $k$ dans $D_1$ (resp. $D_2$) et $c_1$ (resp. $c_2$) le nombre de croix. La caractérisation peut alors s'énoncer ainsi :

\vskip 0.3 cm
\noindent {\bf Caractérisation.}
{\it
Si les $k-1$ premières places sont complètes, la $k$-ième l'est si l'une des conditions suivantes est vérifiée :
\begin{enumerate}
\item en place $k$ dans $D_1$ et $D_2$ il y a une colonne-$y$ et $b_2=b_1+d$ et $c_2=c_1+d'$ (l'une des deux égalités impliquant facilement l'autre) ;
\item en place $k$, il y a une colonne-$x$ croisée dans $D_1$ (i.e. $b_1=0$) et une colonne-$y$ dans $D_2$, et $b_2=d$ ;
\item en place $k$, il y a une colonne-$y$ dans $D_1$ et une colonne-$x$ blanche dans $D_2$ ($c_2=0$), et $c_1=-d'$.
\end{enumerate}}

Pour des raisons de symétrie, nous nous sommes restreint au cas o\`u nous avons une colonne-$y$ dans $T_1+S_2$. Cette caractérisation se prouve aisément en étudiant chacun des cas pouvant se produire.

\subsection{Application}

Une fois ce cadre posé, nous pouvons avancer dans la preuve du Théorème 3. Nous voulons montrer qu'un dessin $D$ est différent de tous ses descendants. Nous pouvons nous restreindre aux descendants qui ont la même forme que lui, i.e. les colonnes-$y$ aux mêmes places. Soit $D'$ un tel descendant. Nous démontrons alors  les lemmes suivants :

\vskip 0.2 cm
\noindent {\bf Lemme 2.}
{\it Si nous avons complétude sur les $k$ premières places le long de la chaîne entre $D$ et $D'$, alors la somme des $d$ le long de la chaîne vaut zéro, de même que la somme des $d'$ ($d$ et $d'$ relatifs aux $k$ premières places).}

\vskip 0.2 cm
Ce lemme nous permet alors de démontrer le :

\vskip 0.2 cm
\noindent {\bf Lemme 3.}
{\it Si nous avons complétude sur les $k$ premières places entre $D$ et $D'$, alors ces deux dessins sont identiques sur ces $k$ places.}

\vskip 0.2 cm
On en déduit alors ce dernier lemme :

\vskip 0.2 cm
\noindent {\bf Lemme 4.} 
{\it S'il n'y a pas complétude totale entre $D$ et $D'$, alors $D\neq D'$, ce qui implique le Théorème 3.}

\vskip 0.2 cm
Il suffit donc de montrer qu'il ne peut y avoir totale complétude entre $D$ et $D'$. C'est l'objet du paragraphe suivant.

\subsection{Fin de la preuve}

Nous montrons en fait qu'à chaque génération, il n'y a pas complétude.

Notons encore $D_1=(S_1,T_1)$ et $D_2=(S_2,T_2)$ deux dessins, père et fils.

Si $D_1$ et $D_2$ ont la même forme, le résultat est évident. Nous supposons donc que $D_1$ et $D_2$ sont de formes différentes et raisonnons par l'absurde en supposant qu'il y a complétude. En regardant la place la plus à gauche o\`u la forme change, on se ramène à un changement de forme en place 1. Les seuls cas o\`u la non-complétude n'est pas évidente sont les cas suivants (remarquer qu'ici $d=d'=0$) :
\vskip 0.5 cm

\centerline{
\epsffile{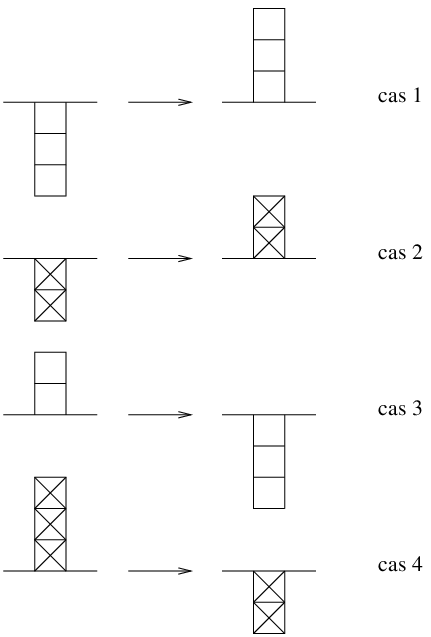}}
 
\vskip 0.2 cm
\noindent {\it Remarque 2.}
$D_2$ est un fils de $D_1$ si et seulement si flip($D_1$) est un fils de flip($D_2$). Ceci nous permet de nous restreindre aux cas 2 et 4. 

Nous montrons alors que dans ces deux cas il n'y a pas complétude, et ce en utilisant les règles de construction des dessins. Ceci achève la preuve du Théorème 3.

\section{\'Eléments de degré en $x$ nul}

Dans cette partie nous présentons comment les objets introduits dans le cas des équerres permettent de trouver facilement une base monomiale du sous-espace de $M_{\mu}$ de degré en $x$ nul que nous notons $M_{\mu}^0$. Il est démontré dans [3] et [7] que la dimension de cet espace est $n!/\mu!$. La base que nous obtenons est liée à une famille introduite dans [3] mais ici non seulement la construction est directe mais en plus on applique les dérivées monomiales à $\Delta_{\mu}$ lui-même.

Soit $\mu$ une partition quelconque de $n$. Ici un dessin est constitué de $n-1$ barres. L'ensemble des couples (profondeur,hauteur) des barres correspondent à l'ensemble des biexposants de $\mu$. Les règles pour placer les barres et les croix dans celles-ci sont les suivantes :
\begin{enumerate}
\item les barres ayant le même nombre de cases-$x$ sont rangées par hauteur décrois\-sante ;
\item il y a des croix dans chaque case-$x$ ;
\item si une barre est placée à gauche d'une barre ayant plus de cases-$x$ et $q$ cases-$y$, alors la première doit avoir au moins $q+1$ cases-$y$ blanches.
\end{enumerate}  

Voici un exemple de dessin :
\vskip 0.5 cm

\centerline{
\epsffile{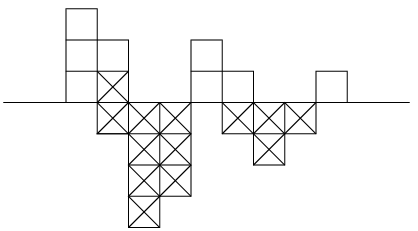}}
\vskip 0.3 cm

\noindent
associé à la partition :
\vskip 0.5 cm

\centerline{
\epsffile{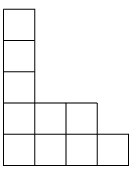}}

\vskip 0.3 cm

\`A un dessin $D$ nous associons $S$  (resp. $T$) l'ensemble des cases croisées (resp. blanches). Nous notons $\partial_S$ l'opérateur de dérivation associé à $S$ et $M_T$ le monôme associé à $T$ (pour chaque case-$y$ en place $i$, nous comptons un $x_i$ et nous faisons le produit sur toutes ses cases). Nous notons enfin ${\cal S}$ l'ensemble des $S$ ainsi construits.  

\vskip 0.3 cm
\noindent {\bf Théorème 4.}
{\it La famille $\{\partial_S\Delta_{\mu}\}_{S\in{\cal S}}$ est une base pour $M_{\mu}^0$.}

Pour cela, nous vérifions que le cardinal de ${\cal S}$ est $n!/\mu!$, nous prouvons que le monôme dominant de $\partial_S\Delta_{\mu}$ pour l'ordre lexicographique est $M_T$ et nous concluons en vérifiant que tous les $T$ sont distincts.

\vskip 1 cm
\noindent {\bf \large Références}
\small
\vskip 0.5 cm
\begin{enumerate}

\item E. Allen, {\it The decomposition of a bigraded left regular representation of the diagonal action of $S_n$}, J. Comb. Theory A, {\bf 71} (1995), 97-111.

\item J.-C. Aval, {\it Monomial bases related to the $n!$ conjecture}, soumis, 27 pages. 

\item N. Bergeron and A. M. Garsia, {\it On certain spaces of harmonic polynomials}, Con\-temporary Mathematics, {\bf 138} (1992), 51-86.

\item Louis Comtet, {\it Analyse Combinatoire}, Presses Universitaires de France, Paris, 1970.

\item A. M. Garsia and M. Haiman, {\it Orbit harmonics and graded representation}, in ``Laboratoire de combinatoire et d'informatique math\'ematique, UQAM collection'' (S. Brlek, Ed), à paraître.

\item A. M. Garsia and M. Haiman, {\it A graded representation model for Macdonald's polynomials}, Proc. Natl. Acad. Sci., {\bf 90} (1993), 3607-3610.

\item A. M. Garsia and M. Haiman, {\it Some natural bigraded $S_n$-modules and $q,t$-Kostka coefficients}, Elec. J. of Comb. 3 (no. 2) (1996), R24. 

\item A. M. Garsia and J. Remmel, {\it Plethistic formulas and positivity of for $q,t$-Kostka polynomials}, In Mathematical Essays in Honor of Gian-Carlo Rota (Cambridge, MA, 1996), Birkh\"auser Boston, Boston, MA (1998), 245-262.

\item A. M. Garsia and G. Tesler, {\it Plethistic formulas for Macdonald $q,t$-Kostka coefficients}, Advances in Math., {\bf 123} (1996), 144-222.

\item M. Haiman, {\it Macdonald polynomials and geometry}, preprint.

\item A. N. Kirillov and M. Noumi, {\it Affine Hecke algebras and raising operators for Macdonald polynomials}, Duke Math. J., {\bf 93} (1998), 1-39.

\item F. Knop, {\it Integrality of two variable Kostka functions}, J. Reine Angew. Math., {\bf 482} (1997), 177-189.

\item I. G. Macdonald, {\it A new class of symmetric functions}, Actes du $20^e$ S\'eminaire Lotharingien, Publ. I.R.M.A. Strasbourg (1988), 131-171.

\item E. Reiner, {\it A Proof of the $n!$ Conjecture for Generalized Hooks}, J. Comb. Theory A, {\bf 75} (1996), 1-22.

\item S. Sahi, {\it Interpolation, integrality, and a generalization of Macdonald's polynomials}, Internat. Math. Res Notices, {\bf 10} (1996), 457-471.

\end{enumerate}

\end{document}